\documentstyle[12pt]{article}

\input psfig.sty
\newtheorem{df}{ \sc Definition}[section]
\newtheorem{as}[df]{ \sc Assumption}
\newtheorem{ex}[df]{ \it Example}
\newtheorem{pr}[df]{ \sc Proposition}

\newtheorem{th}[df]{ \sc Theorem}
\newtheorem{cor}[df]{ \sc Corollary}
\newtheorem{re}[df]{ \it Remark}

\def\nat{\smash\Join\smash}
\def\n{{\oplus n}}
\def\codim{{\rm codim}}
\def\mpr#1{\;\smash{\mathop{\hbox to 20pt{\rightarrowfill}}\limits^{#1}}\;}
\def\epi#1{\;\smash{\mathop{\hbox to 20pt{\rightarrowfill}\hskip
-13pt\rightarrow}\limits^{#1\,}}\;}
\def\epii{\smash{\mathop{\hbox to 14pt{\rightarrowfill}\hskip
-13pt\rightarrow}}}

\def\mpl#1{\;\smash{\mathop{\hbox to 20pt{\leftarrowfill}}\limits^{#1}}\;}

\def\B{{ B}}
\def\C{{\cal PR}}
\def\Ct{{\bf C}}

\def\Zt{{\bf Z}}
\def\Nt{{\bf N}}

\def\Hi{H^\flat}
\def\Ho{H^\sharp}
\def\Pii{P^\flat}
\def\Po{P^\sharp}
\def\Proof{\noindent{\it Proof. }}

\begin{document}

\author{Bogdan Bojarski, Andrzej Weber
\footnote{Supported by KBN grant  1 P03A 005 26
 }}
\title{CORRESPONDENCES AND INDEX }
\maketitle
\begin{abstract} We define  certain class of correspondences of
polarized representations of $C^*$-algebras. Our correspondences
are modeled on the spaces of boundary values of elliptic operators
on bordisms joining two manifolds. In this setup we define the
index. The main subject of the paper is the additivity of the
index.\end{abstract}

\section{Introduction}

Let $X$ be a closed manifold. Suppose it is decomposed into a sum
of two manifolds $X_+$, $X_-$ glued along the common boundary
$$\partial X_+=\partial X_-=M\,.$$ Let $$ D :C^\infty(X;\xi
)\rightarrow C^\infty(X;\eta)$$ be an elliptic operator of the
first order. We assume that it possesses the unique extension
property: if $Df=0$ and $f_{|M}=0$ then $f=0$. In what follows we
will consider only elliptic operators of the first order such that
$D$ and $D^*$ have the unique extension property.

One defines the spaces $H_\epsilon( D )\subset L^2(M;\xi)$ for
$\epsilon\in\{+,-\}$, which are the closures of the spaces of
boundary values of solutions of $Df=0$ on the manifolds
$X_\epsilon$ with boundary $\partial X_\epsilon=M$. The space
$H_\epsilon(D)$ is defined to be the closure of : $$\{f\in
C^\infty(M;\xi)\,:\;\exists \tilde f\in
C^\infty(X_\epsilon;\xi),\; f=\tilde f_{|M},\; \; D (\tilde
f)=0\;\}$$ in $ L^2(M;\xi)$. The pair of spaces $H_\pm( D )$ is a
Fredholm pair, \cite{Bo1}. There are associated Calder\'on
projectors $P_+(D)$ and $P_-(D)$, see \cite{Sel}.

To organize somehow the set of possible Cauchy data we will
introduce certain algebraic object. We fix a $C^*$-algebra $\B$,
which is the algebra of functions on $M$ in our case. Suppose it
acts on a Hilbert space $H$. Now we consider Fredholm pairs in
$H$. In our case $H=L^2(M;\xi)$ and one of the possible Fredholm
pairs is $H_\pm(D)$. Note that this pair is not arbitrary. It has
a property which we called {\it good}. A Fredholm pair is good if
(roughly speaking) it remains to be Fredholm after conjugation
with functions, see \S\ref{Good}. These pairs act naturally on
$K_{1}(M)$. Nevertheless the concept of a good Fredholm pair is
not convenient to manipulate, thus we restrict our attention to
the pairs of geometric origin, see \S\ref{Adm}. We call them {\it
admissible}. They are the pairs of subspaces which are images of
projectors which almost commute with the actions of the algebra
$\B$. This concept allows to extract the relevant
analytico--functional information out of the Cauchy data. Further
a Morse decomposition of a manifold is translated into this
language.

Our paper is devoted to the study of the cut and paste technique
on manifolds and its effect on indices. The spirit of these
constructions comes from the earlier papers \cite{Bo1}--\cite{Bo3}
or \cite{BW1}. According to the topological and conformal field
theory we investigate the behaviour of the index of a differential
operator on a manifold composed from bordisms
$$X=X_0\cup_{M_1}X_1\cup_{M_2}\dots\cup_{M_{m-1}}X_{m-1}\cup_{M_m}X_m\,.$$
We think of $M_i$'s as objects and we treat bordisms of manifolds
as morphisms. Starting from this geometric background we introduce
a category $\C$, whose objects are {\it polarized
representations}. The algebra $\B$ may vary. We keep in mind that
such objects arise when:
 \begin{itemize}
 \item $\B$ is an algebra of functions on a manifold $M$,
 \item there is given a vector bundle $\xi$ over $M$, then
$H=L^2(M;\xi)$ is a representation of $\B$,
 \item there is given a pseudodifferential projector in $H$.
\end{itemize}
 The morphisms in $\C$ are certain correspondences,
i.e.~linear subspaces in the product of the source and the target.
A particular case of principal value for our theory are the
correspondences coming from bordisms of manifolds equipped with an
elliptic operator. Precisely: suppose we are given a manifold $W$
with a boundary $\partial W=M_1\sqcup M_2$. Moreover, suppose that
there is given an elliptic operator of the first order acting on
the sections of a vector bundle $\xi$ over $W$. Then the space of
the boundary values of the Cauchy data of solutions is a linear
subspace in $L^2(M_1;\xi_{|M_1})\oplus L^2(M_2;\xi_{|M_2})$. In
another words it is a correspondence from $L^2(M_1;\xi_{|M_1})$ to
$L^2(M_2;\xi_{|M_2})$.

{\it Basic example:} The following example is instructive and
serves as the model situation (see \cite{BWe}): Let $W=\{z\in
\Ct\,:\, r_1\geq |z| \geq r_2\}$ be a ring domain and let $D$ be
the Cauchy-Riemann operator. The space $L^2(M_i)$ for $i=1,2$ is
identified with the space of sequences $\{a_n\}_{n\in \Zt}$, such
that $\sum_{n\in \Zt} |a_n|^2r_i^{2n}<\infty$. The sequence
$\{a_n\}$ defines the function on $M_i$ given by the formula
$f(z)=\sum_{n\in\Zt}a_nz^n$. The subspace of the boundary values
of holomorphic functions on $W$ is identified with
 $$\Big\{(\{a_n\},\{b_n\}) \,:\, \Sigma_{n\in
\Zt} |a_n|^2r_1^{2n}<\infty\,,\;\Sigma_{n\in \Zt}
|b_n|^2r_2^{2n}<\infty\; {\rm and}\;a_n=b_n\Big\}\,.$$
 It can be
treated as the graph of an unbounded operator $\Phi:L^2(M_1) \to
L^2(M_2)$. When we restrict $\Phi$ to the space $L^2(M_1)^\sharp$
consisting of the functions with coefficients $a_n=0$ for $n<0$ we
obtain a compact operator. On the other hand the inverse operator
$\Phi^{-1}:L^2(M_2) \to L^2(M_1)$ is compact when restricted to
$L^2(M_2)^\flat$, the space consisting of the functions with
coefficients $a_n=0$ for $n\geq 0$.

 The Riemann-Hilbert transmission problem of the Cauchy data
across a hypersurface is a model for another class of morphisms.
These are called {\it twists}. Our approach allows us to treat
bordisms and twists in a uniform way. We calculate the global
index of an elliptic operator in terms of local indices depending
only on the pieces of the decomposed manifold (see Theorems
\ref{chain-ind} and \ref{add+twist}). An interesting phenomenon
occurs. The index is not additive with respect to the composition
of bordisms. Instead each composition creates a contribution to
the global index (Theorem \ref{strange}):
 $$L_1,\,L_2\leadsto L_2\circ L_1 +\delta(L_1,L_2)\,.$$
 In the
geometric situation this contribution might be nonzero for example
when a closed manifold is created as an effect of composition of
bordisms. One can show that if the bordisms in $\C$ come from
connected geometric bordisms supporting elliptic operators with
the unique extension property then the index is additive. The
contributions coming from twists are equivalent to the effects of
pairings in the odd $K$-theory, Theorem \ref{chain-zap}.

It's a good moment now to expose a fundamental role of the
splitting of the Hilbert space into a direct sum. The need of
introducing a splitting was clear already in \cite{Bo1}:
\begin{itemize}
\item It was used to the study of Fredholm pairs with
application to the Riemann-Hilbert transmission problem in
\cite{Bo1}
\item Splitting also came into light
in the paper of Kasparov \cite{Ka}, who introduced a homological
$K$-theory built from the Hilbert modules. The program of
noncommutative geometry of A.Connes develops this idea, \cite{Co1,
Co2} .
\item Splitting plays an important role in the theory of loop
groups in \cite{PSe}.
\item There is also a number of papers in which surgery of the
Dirac operator is studied. Splitting serves as a boundary
condition, see e.g.~\cite{DZ}, \cite{SW}. These papers originate
from \cite{APS}.
\end{itemize}

In the present paper we omit the technicalities and problems
arising for a general elliptic operator. We concentrate on the
purely functional calculus of correspondences. This is mainly the
linear algebra.

\section{\label{p1}Fredholm pairs}

Let us first summarize some facts about Fredholm pairs. We will
follow \cite{Bo1}-\cite{Bo3}. Suppose that $H_+$ and $H_-$ are two
closed subspaces of a Hilbert space, such that $H_++H_-$ is also
closed and
\begin{itemize}
\item $H_+\cap H_-$ is of finite dimension,
\item $H_++ H_-$ is of finite codimension.\end{itemize}
We assume that both spaces have infinite dimension. Then we say
that the pair $(H_+,H_-)=H_\pm$ is Fredholm. We define its index
$$Ind(H_\pm)=\dim(H_+\cap H_-)-\codim(H_++ H_-)\,.$$
The following statements follow from easy linear algebra.

\begin{pr}\label{kryt1} A pair $H_\pm$ is Fredholm if, and only if the map
$$\iota:H_+\oplus H_-\rightarrow H$$
induced by the inclusions is a Fredholm operator. Moreover the
indices are equal: $$Ind(H_\pm)=ind(\iota)\,.$$\end{pr}

Here $Ind$ denotes the index of a pair, whereas $ind$ stands for
the index of an operator. Suppose that $H$ is decomposed into a
direct sum
$$H=\Hi\oplus \Ho\,.$$
We may assume that this decomposition is given by a symmetry
$S$:~a "sign" or "signature" operator. Let $\Pii$ and $\Po$ be the
corresponding projectors. We can write $S=\Po-\Pii$. We easily
have:

\begin{pr}\label{kryt2} If $H_\pm$ is a pair with $H_+=\Ho$, then it
is Fredholm if and only if the restriction
$\Pii_{|H_-}:H_-\rightarrow \Hi$ is a Fredholm operator. Moreover
the indices are equal: $$Ind(H_\pm)=ind(\Pii_{|H_-})\,.$$\end{pr}

Let ${\cal I}\subset L(H)$ be an ideal which is lies between the
ideal of finite rank operators and the ideal of compact operators
$${\cal F}\subset{\cal I}\subset{\cal K}\,.$$
Define $GL(\Pii,{\cal I})\subset GL(H)$ to be the set of the
invertible automorphisms of $H$ commuting with $\Pii$ up to the
ideal $\cal I$. We will say that $\phi$ almost commutes with
$\Pii$ or we will write $\phi\Pii\sim \Pii \phi$. Obviously
$GL(\Pii,{\cal I})=GL(\Po,{\cal I})=GL(S,{\cal I})$.
 We
have the following description of Fredholm pairs stated in
\cite{Bo1}. (The proof is again an easy linear algebra.)

\begin{th} \label{bo} Let $H_\pm$ be a Fredholm pair
with $H_+=\Ho$. Then there exists a complement $\Hi$ (that is
$\Hi\oplus \Ho=H$) and there exists $\phi\in GL(\Pii,{\cal I})$,
such that $H_-=\phi(\Hi)$. If $H_\pm$ is given by a pair of
projectors $P_\pm$ satisfying $P_-+P_--1\in {\cal I}$, then we can
take $\Hi=ker\,P_+$. Moreover, the operator $\phi\Pii+\Po$ is
Fredholm and
$$ind(\phi\Pii+\Po)=Ind(H_\pm)\,.$$
The map
$$\widetilde{ind}:GL(\Pii,{\cal I})\rightarrow \Zt$$
$$\widetilde{ind}(\phi)=ind(\phi\Pii+\Po)$$
is a group homomorphism.
\end{th}
It follows that
$$ind(\phi\Pii+\Po)= ind(\Pii\phi:\Hi\rightarrow
\Hi)=ind(\Po\phi^{-1}:\Ho\rightarrow \Ho)\,.$$

\section{\label{egz}Index formula for a decomposed manifold}

The main example of a Fredholm pair is the following. Let $D$ be
an elliptic operator on $X=X_+\cup_M X_-$. Then the pair of
boundary value spaces $H_\pm(D)$ (as defined in the introduction)
is a Fredholm pair.

\begin{as} [Unique Extension Property] \label{uep} Let $\epsilon=+$ or $-$ and
let $f\in
C^\infty(X_\epsilon;\xi)$.
 If $Df=0$ and $f_{|M}=0$ then $f=0$.
\end{as}

If $D$ has the unique extension property, then $${\rm
ker}(D)\simeq H_+(D)\cap H_-(D)\,.$$ This formula is easy to
explain: a global solution restricted to $M$ lies in $H_+(D)\cap
H_-(D)$. On the other hand if a section $f$ of $\xi$ over $M$ can
be extended to both $X_+$ and $X_-$, such that the extensions are
solutions of $Df=0$ then we can glue them to obtain a global
solution. The unique extension property is necessary, because we
need to know that a solution is determined by its restriction to
$M$. Following the reasoning in \cite{Bo1}, with the Assumption
\ref{uep} for $D$ and $D^*$ we have:

\begin{cor}\label{both}
$$Ind(H_\pm(D))=ind(D)\,.$$ \end{cor}
For a rigorous proof see \cite{BW2}, \S24 for Dirac type
operators.

\begin{re} \rm It may happen that $D$ does not have the unique extension
property. This is so for example when $X$ is not connected. Then
the Cauchy data $H_\pm(D)$ do not say anything about the index of
the operator $D$ on the components of $X$ disjoined with $M$.
There are also known elliptic operators without the unique
extension property on connected manifolds, \cite{P}, \cite{A}. It
is difficult to characterize the class of all operators $D$ with
the unique extension property. Nevertheless the most relevant are
Cauchy-Riemann  and Dirac type operators. These operators have the
unique extension property.

\end{re}

\section{\label{Good}Good Fredholm pairs}

Suppose there is given an algebra $\B$ and its representation
$\rho$ in a Hilbert space $H$. For a Fredholm pair $H_\pm$ in $H$
and an invertible matrix $A\in GL_n(\B)$ we define a new pair of
subspaces $A\nat H_\pm$ in $H^\n $. We set
 $$(A\nat H_\pm)_-=\rho A(H_-^\n )\qquad (A\nat H_\pm)_+=H_+^\n \,.$$
(As usually we treat $\rho A$ as an automorphism of $H^\n $.)

\begin{df} \rm Let $\B$ be a $C^*$-algebra which acts on a Hilbert space $H$.
A {\it good Fredholm pair} is a pair of subspaces $(H_+,H_-)$ in
$H$, such that for any invertible matrix $A\in GL(n;\B)$ the pair
$A\nat H_\pm$ is a Fredholm pair.
\end{df}

We will see that the pair of boundary values $H_\pm(D)\subset
H=L^2(M;\xi)$ for the operator $D$ considered in the introduction
is good.

\begin{ex}\label{RH}\rm [Main example: Riemann-Hilbert problem]
Consider the following problem: there is given a matrix-valued
function $A:M\to GL_n(\Ct)$. We look for the sequence
$(s_\pm^1,\dots,s_\pm^n)$ of solutions of $Ds=0$ on $X_\pm$
satisfying the transmission condition on $M$
 $$A(s_-^1,\dots,s_-^n)=(s_+^1,\dots,s_+^n)\,.$$
A Fredholm operator is related to this problem and we study its
index, see \S\ref{weird}. On the other hand the matrix $A$ treated
as the gluing data defines an $n$-dimensional vector bundle
$\Theta_X^A$ over $X$. Then
 $$Ind(A\nat H_\pm( D ))=ind( D \otimes\Theta_X^A)\,.$$
This formula was obtained in \cite{BW1}, \S1 under the assumption
that $D$ has a product form along $M$.
\end{ex}

\begin{cor} \label{Hgfp} For the elliptic operator $ D $
the pair $H_\pm( D )\subset L^2(M;\xi)$ is a good Fredholm
pair.\end{cor}

\begin{re}\rm Consider the differential in the Mayer-Vietoris
exact sequence of $X=X_+\cup_M X_-$
 $$\delta:K_0(X)\to K_{-1}(M)\,.$$
The operator $D$ defines a class $[D]\in K_0(X)$. The element
$\delta[D]$ can be recovered from the good Fredholm pair
$H_\pm(D)\subset L^2(M;\xi)$. Note that the pair $H_\pm( D )$
encodes more information. One can recover the index of the
original operator. We describe the map $\delta$ via duality,
therefore we neglect the torsion of $K$-theory. The construction
is the following: for an element $a\in K^1(M)$ we define the value
of the pairing
 $$\langle \delta[D],a\rangle=\langle [D],\partial a\rangle\,.$$
The element $a$ is represented by a matrix $A\in
GL_n(C^\infty(M))$. Then
 $$\langle [D],\partial a\rangle
=ind(D\otimes\Theta_X^A)-n\,ind(D)\,,$$
 where $\Theta_X^A$ is the bundle defined in Example \ref{RH}.
Now
 $$\langle [D],\partial a\rangle=Ind(A\nat H_\pm( D ))-n Ind(H_\pm( D
 ))\,.$$

\end{re}

\section{\label{Adm}Admissible Fredholm pairs}

The following can be related to the paper of Birman and Solomyak
\cite{BS} who introduced the name {\it admissible} for the
subspaces which are the images of pseudodifferential projectors.
Suppose that $\xi$ is a vector bundle over a manifold $M$. We
consider Fredholm pairs $H_\pm$ in $H=L^2(M;\xi)$ such that the
subspaces $H_\pm$ are images of pseudodifferential projectors
$P_\pm$ with symbols satisfying
 $$\sigma(P_+)+\sigma(P_-)=1\,.$$
We would like to free ourselves from the geometric context and
state admissibility condition in an abstract way. We assume that
$H$ is an abstract Hilbert space with a representation of an
algebra $\B$, which is the algebra of functions on $M$ in the
geometric case. The condition that $P_\pm$ is pseudodifferential
we substitute by the condition: $P_\pm$ commutes with the algebra
action up to compact operators. We are ready now to give a
definition:

\begin{df} \rm We say that a pair of subspaces $H_\pm$ is an {\it admissible
Fredholm pair} if there exist a pair of projectors $P_\epsilon$
for $\epsilon\in\{+,-\}$, such that $H_\epsilon=im\,P_\epsilon$
and $P_\epsilon$ commutes with the action of $\B$ up to compact
operators. Moreover, we assume that $P_++P_--1$ is a compact
operator.
\end{df}

\begin{pr} \label{eachad} Each admissible Fredholm pair is a good Fredholm pair.\end{pr}

\Proof Set $K=P_++P_--1$. If $v\in H_+\cap H_-$, then $K(v)=v$.
Since $K$ is a compact operator, $\dim(H_+\cap H_-)<\infty$. To
prove that $H_++ H_-$ is closed and of finite codimension, note
that $im(P_++ P_-)\subset H_++ H_-$. Since $P_++ P_-$ is Fredholm
its image is closed and of finite codimension. This way we have
shown that $H_\pm$ is a Fredholm pair. Now, if we conjugate
$P_+^\n $ by $\rho A$ we obtain again an almost complementary pair
of projectors. Thus $A\nat H_\pm$ is a Fredholm pair as
well.\hfill$\Box$ \vskip10pt

We denote by $AFP(\B)$ the set of good Fredholm pairs
divided by the equivalence relation generated by homotopies and stabilization
with respect to the direct sum. We also consider as trivial the pairs
associated to projectors strictly satisfying $P_++P_-=1$ and
commuting with the action of $\B$. In another words these are just
direct sums of two representations of $\B$. It is not hard to show
that
\begin{pr}
$$AFP(\B)\simeq K^1(\B)\oplus \Zt\,.$$\end{pr}
\Proof We have the following natural transformation:
$$\matrix{\beta:&AFP(M)&\rightarrow & K_{1}(M)\cr
&(H,P_\pm)&\mapsto & (H,S_+)\cr}\,.$$ Here $S_+=2P_+-1$ is s just
the symmetry defined by $P_+$. We remind that the objects
generating $ K_{1}(M)$ are odd Fredholm modules, see \cite{Co2},
pp 287-289. This procedure is simply forgetting about $P_-$. We
can recover $P_-$  (up to homotopy) by fixing the index of the
pair, i.e $\beta\oplus Ind$ is the isomorphism we are looking for.
Precisely, the pseudodifferential projector is determined up to
homotopy by its symbol and the index, see \cite{BW2}.
\hfill$\Box$\vskip10pt

\section{\label{ssp}Splittings and polarization}
We adopt the concepts of splitting and polarization to our
situation.

\begin{df}\rm Let $H$ be a representation of a $\Ct^*$-algebra $\B$ in a
Hilbert space. A {\it splitting} of $H$ is a decomposition
$$H=\Hi\oplus\Ho\,,$$
such that the projectors on the subspaces $\Pii$, $\Po$ commute
with the action of $\B$ up to compact operators.
\end{df}

The basic example of a splitting is the one coming from a
pseudodifferential projector. Another equivalent way of defining a
splitting(as in \cite{Bo2}) is to distinguish a symmetry $S$,
almost commuting with the action of $\B$. Then $\Hi$ is the
eigenspace of $-1$ and $\Ho$ is the eigenspace of $1$. Then we may
think of $H$ as a superspace, but we have to remember that the
action of $\B$ does not preserve the grading.

\begin{df}\rm
In the set of splittings we introduce an equivalence relation: two
splittings are equivalent if the corresponding projectors coincide
up to compact operators. An equivalence class of the above
relation is called a {\it polarization} of $H$.\end{df}

Informally we can say, that polarization is a generalization of
the symbol of a pseudodifferential projector.

\begin{ex}\label{sym->pol}\rm Let $\xi\rightarrow M$ be a complex vector
bundle over a manifold. Let $\widetilde\xi$ be the pull back of
$\xi$ to $T^*M\setminus\{0\}$. Suppose
$p:\widetilde\xi\rightarrow\widetilde\xi$ is a bundle map which is
a projector (hence $p$ is homogeneous of degree 0) . Then $p$
defines a polarization of $L^2(M;\xi)$. Just take a
pseudodifferential projector $P=\Po$ with $\sigma(P)=p$ and set
$$\Hi=ker\,P\,,\qquad\Ho=im\,P\,.$$\end{ex}

\begin{ex}\label{pair->pol}\rm Suppose $(H_+,H_-)$ is an admissible Fredholm
pair given by projectors $(P_+,P_-)$. Then the polarizations
associated with $P_+$ and $1-P_-$ coincide. This way an admissible
Fredholm pair defines a polarization. Furthermore each
polarization defines an element of $K_1(\B)$.\end{ex}

Intuitively polarizations can be treated as a kind of orientations
dividing $H$ into the upper half and lower half. Such a tool was
used in \cite{DZ} to split the index of a family of Dirac
operators. (In \cite{DZ} splittings were called generalized
spectral sections.) Polarizations were discussed in the lectures
of G.~Segal (see \cite{Seg}, Lecture 2).

\section{\label{cbt}Correspondences, bordisms, twists}
\begin{df}\label{cat}\rm We consider the category $\C$ having the following
 objects and morphisms
\begin{itemize}
\item $Ob(\C)\;=\;$Hilbert spaces (possibly of
finite dimension) with a representation of some $\Ct^*$-algebra $\B$
and with a distinguished polarization,
\item $Mor_{\C}(H_1,H_2)\;=\;$ closed linear subspaces $L\subset H_1\oplus
H_2$, such that the pair $(L,\Hi_1\oplus\Ho_2)$ is Fredholm.
\end{itemize} We write also $H_1\mpr{L} H_2$.
\end{df}
In particular
$$Mor_{\C}(H,0)\subset Grass(H)\supset Mor_{\C}(0,H)\,.$$
By Proposition \ref{kryt2} a subspace $L\subset H_1\oplus H_2$ is a
morphism if and only if
$$\Pi=\Po_1\oplus\Pii_2:L\rightarrow \Ho_1\oplus\Hi_2$$
is a Fredholm operator.
The composition in $\C$ is the standard composition of correspondences:
$$L_{1}\subset H_1\oplus H_2\,,\qquad L_{2}\subset H_2\oplus H_3\,,$$
$$L_{2}\circ L_{1}=\{(x,z)\in H_1\oplus H_3\;:\;\exists y\in H_2\,,\;(x,y)\in
L_{1}\,,\;(y,z)\in L_{2}\;\}\,.$$ In another words the morphisms
are certain correspondences or relations, as they were called in
\cite{Bo1}. Our approach also fits to the ideas of the topological
field theory as presented in \cite{Seg}.

\begin{pr} \label{comor} The composition of morphism is a morphism.\end{pr}

\Proof Let $L_{1}\in Mor_\C(H_1,H_2)$ and $L_{2}\in
Mor_\C(H_2,H_3)$. A simple linear algebra argument shows that
\begin{itemize}
\item
the kernel of
$$\Pi_{13}:L_{2}\circ L_{1}\rightarrow \Ho_1\oplus\Hi_3$$
is a quotient of $ker(\Pi_{12})\oplus ker(\Pi_{23})$,
\item the cokernel of $\Pi_{13}$ is a subspace of
$coker(\Pi_{23})\oplus coker(\Pi_{12})$.\hfill$\Box$
\end{itemize}

The role of polarizations in the definition of morphisms is clear
and the algebra actions are involved implicitly. In fact, the
object which plays the crucial role is the algebra of operators
commuting with $\Po$ up to compact operators, i.e.~the odd
universal algebra. The role of this algebra was emphasized in
$\cite{Bo2}$. However, in the further presentation we prefer to
expose the geometric origin of our construction and keep the name
$B$.

We have two special classes of morphisms in $\C$:

\begin{df} \label{twist}\rm A subspace $L\subset H\oplus H$ is a
{\it twist} if it is the graph of a linear isomorphism $\phi\in
GL(\Po,{\cal K})\subset GL(H)$ commuting with the polarization
projectors up to compact operators.\end{df}

\begin{pr} For a twist $L=graph(\phi)\subset H\oplus
H$ the pair $(L,\Hi\oplus\Ho)$ is Fredholm, i.e. $L\in Mor_\C(H,H)$.
\end{pr}

\Proof To show that $(L,\Hi\oplus\Ho)$ is a Fredholm pair let us show
that the projection
$$\Pi=\Po\oplus\Pii:L\rightarrow \Ho\oplus\Hi\subset H\oplus H$$
is a Fredholm operator. Indeed, $L$ is parameterized by
$$(1,\phi):H\rightarrow L\subset H\oplus H\,.$$ The composition of these maps
is
equal to
$$F=\Po\oplus\Pii\phi\,.$$
Since $\phi$ almost commutes with $\Pii$ the map $F$ has a
parametrix $\widetilde F=\Po\oplus\Pii\phi^{-1}\,.$\hfill$\Box$

\begin{df} \label{bordism} \rm A subspace $L\subset H_1\oplus H_2$
is a bordism if $L$ is the image of a projector $P_L$, such that
$$P_L\sim \Po_1\oplus \Pii_2\,.$$
\end{df}

By \ref{eachad} for any $P_L\sim\Po_1\oplus \Pii_2$ the pair
$(L,\Hi_1\oplus\Ho_2)$ is Fredholm. The motivation for the
Definition \ref{bordism} is the following:

\begin{ex} \label{bor-geom}\rm Let $X$ be a bordism between closed
manifolds $M_1$ and $M_2$, i.e. $$\partial X=M_1 \sqcup M_2\,.$$
Suppose that $D:C^\infty(X;\xi)\rightarrow C^\infty(X;\eta)$ is an
elliptic operator of the first order. Then the symbols of
Calder\'on projectors define polarizations of $H_1=L^2(M_1;\xi)$
and $H_2=L^2(M_2;\xi)$, see Example \ref{sym->pol}. We reverse the
polarization on $M_2$, i.e.~we switch the roles of $\Hi$ and
$\Ho$. Let $L\subset L^2(M_1;\xi)\oplus L^2(M_2;\xi)$ be the
closure of the space of boundary values of solutions of $Ds=0$.
Then $L\in Mor_\C(H_1,H_2)$ is a bordism in $\C$. This procedure
indicates the following:
\begin{itemize}
\item the space $L\subset L^2(M_1\sqcup M_2;\xi)=L^2(M_1;\xi)\oplus
L^2(M_2;\xi)$ and the associated Calder\'on projector are {\it
global} objects. One cannot recover them from the separated data
in $L^2(M_1;\xi)$ and $L^2(M_2;\xi)$.
\item but up to compact operators one can {\it localize} the
projector $P_L$ and obtain two projectors acting on $L^2(M_1;\xi)$ and
$L^2(M_2;\xi)$. \end{itemize}
\end{ex}

We note that the following proposition holds:
\begin{pr}\break\begin{enumerate}
\item The composition of bordisms is a bordism.
\item The composition of a bordism and a twist is a bordism.
\item The composition of twists is a twist.\end{enumerate}\end{pr}

\begin{re}\rm Let $H_1\mpr{L_1}H_2\mpr{L_2}H_{3}$ be a
pair of bordisms in $\C$ coming from  geometric bordisms
$$M_1\sim_{X_1}M_2\,,\quad M_2\sim_{X_2}M_3$$
and an elliptic operator on $X_1\cup_{M_2}X_2$, as in Example
\ref{bor-geom}. Then $L_2\circ L_1$ coincides with the space of
the Cauchy data along $\partial(X_1\cup_{M_2}X_2)=M_1\sqcup M_3$
of the solutions of $Ds=0$ on $X_1\cup_{M_2}X_2$.\end{re}

\section{Chains of morphisms}

Now we introduce the notion of a chain. This is a special case of
a Fredholm fan considered in \cite{Bo2} and in \S\ref{k2} below.

A chain of morphisms is a sequence correspondences
$$0\mpr{L_0}H_1\mpr{L_1}H_2\mpr{L_2}\dots\mpr{L_{m-1}}H_m\mpr{L_m}0\,.$$

\begin{ex} \rm Let $(H_+,H_-)$ be an admissible Fredholm pair in $H$.
Then we have a sequence
$$0\mpr{H_-}H\mpr{H_+}0$$
which is a chain of bordisms with respect to the polarization
defined by $\Po=P_+$ (or $1-P_-$), see Example
\ref{pair->pol}.\end{ex}

\begin{ex} \rm Each morphism in $L\in Mor_\C(H_1,H_2)$ can be
completed to a chain
$$0\mpr{L_1}H_1\mpr L H_2 \mpr{L_2} 0\,.$$
Just take $L_1=(0\oplus\Hi_1)\subset (0\oplus H_1)$ and
$L_2=(\Ho_2\oplus 0)\subset (H_2\oplus 0)$.\end{ex}

\begin{ex} \label{geom-chain} \rm It is proper to explain why we are interested in
chains of morphisms. Suppose there is given a closed manifold
which is composed of usual bordisms
$$X=X_0\cup_{M_1}X_1\cup_{M_2}\dots\cup_{M_{m-1}}X_{m-1}\cup_{M_m}X_m\,.$$
We treat the manifolds $M_i$ as objects and  bordisms
$$M_{i-1}\sim_{X_i}M_i$$ as morphisms. In particular
$$\emptyset\sim_{X_1}M_1\qquad{\rm and }\qquad M_m\sim_{X_m}\emptyset\,.$$
Let $D:C^\infty(X;\xi)\rightarrow C^\infty(X;\eta)$ be an elliptic
operator of the first order.
This geometric situation gives rise to a chain of bordisms in the
category $\C$:
\begin{itemize}
\item $H_i=L^2(M_i;\xi)$ with the action of $\B_i=C(M_i)$ and the
polarization defined by the symbol of Calder\'on projector, as in
\ref{bor-geom},
\item $L_i\subset L^2(M_i;\xi) \oplus L^2(M_{i+1};\xi)$ is the space
of boundary values of the solutions of $Ds=0$ on
$X_i$.\end{itemize}\end{ex}

\section{\label{ioac}Indices in $\C$}
\begin{df}\label{dein} \rm Fix the splittings $S$ of the objects of $\C$.
The pair $(L,\Hi_1\oplus\Ho_2)$ in $H_1\oplus H_2$ is Fredholm by
Definition \ref{cat}. Define the index of a morphism $L\in
Mor_\C(H_1,H_2)$ by the formula:
 $$Ind_{S_1,S_2}(L)\stackrel{\rm def}{=} Ind(L,\Hi_1\oplus\Ho_2)=
ind(\Po_1\oplus\Pii_2:L\rightarrow\Ho_1\oplus \Hi_2)\,.$$ \end{df}

\begin{pr}\label{twin} We have the equality of indices for a twist
\begin{enumerate}
\item $Ind_{S,S}(graph\,\phi)$,
\item index of  $\left(\matrix{1&\Pii\cr\phi&\Po\cr}\right):
H\oplus H\rightarrow H\oplus H$,
\item $\widetilde{ind}(\phi)=ind(\phi\Pii+\Po)=
Ind(\phi(\Hi),\Ho)$ (compare Theorem \ref{bo}),
\end{enumerate}
\end{pr}

\Proof The graph of $\phi$ is parameterized by $(1,\phi)$ and
$\Hi\oplus\Ho$ is parameterized by $(\Pii,\Po)$. Thus by Theorem
\ref{kryt1} the first equality follows. Now we multiply the matrix
(2.) from the left by the symmetry
$\left(\matrix{\Po&\Pii\cr\Pii&\Po\cr}\right)$ and we obtain
$\left(\matrix{\Pii\phi+\Po&0\cr\Po\phi+\Pii&1\cr}\right)\sim
\left(\matrix{\phi\Pii+\Po&0\cr\phi\Po+\Pii&1\cr}\right)$. The
second equality follows.\hfill$\Box$

\begin{re}\label{split-dep} \rm The index of a twist
depends only on the polarization, not on the particular splitting.
This is clear from \ref{twin}.2. It is worthwhile  to point out
that if the twist $\phi=\widetilde A:H^{\oplus n}\to H^{\oplus n}$
is given by a matrix $A\in GL_n(\B)$, then
$$\widetilde{ind}(\widetilde A)=\langle[\widetilde A],[S_{\Hi}]\rangle\,,$$
where $S_{\Hi}$ is the symmetry with respect to $\Hi$ and the
bracket is the pairing in $K$-theory of $K^1(\B)$ with
$K_1(\B)$.\end{re}

On the other hand $Ind_{S_1,S_2}(L)$ does depend on the splitting
for general morphisms.

\begin{re}\rm The index in Example \ref{bor-geom} is equal to
the index of the operator $D$ with the boundary conditions given
by the splittings, as in \cite{APS}.
\end{re}

\begin{re}\rm There are certain morphisms in $\C$ which are
interesting from the point of view of composition. We will say
that $L$ is a {\it special} correspondence if:
\begin{itemize}
\item $L$ is the graph of an injective function $\phi$ defined on
a subspace of $H_1$,
\item the images of the projections of $L$ onto $H_1$ and $H_2$
are dense.
\end{itemize}
(The second condition is equivalent to the first one for the
adjoint correspondence defined as the ortogonal complement
$L^\perp$.) If $L$ is special, then
$$Ind_{S_1,S_2}(L)=Ind(L(\Hi_1),\Ho_2)\,,$$
 where
 $$L(\Hi_1)=\{y\in H_2\,:\;\exists x\in \Hi_1\;\; (x,y)\in L\,\}\,.$$
 Indeed in this case we have
$$L\cap (\Hi_1\oplus \Ho_2)\simeq L(\Hi_1)\cap\Ho_2\quad{\rm and}\quad
L^\perp\cap (\Hi_1\hskip1pt^\perp\oplus
\Ho_2\hskip1pt^\perp)\simeq
L^\perp(\Hi_1\hskip1pt^\perp)\cap\Ho_2\hskip1pt^\perp\,.$$
 Of
course each twist is a special morphism. Another example of a
special morphism is the one which comes from the Cauchy-Riemann
operator. In general, we obtain a special morphism if the operator
(and its adjoint) satisfies the following:
\begin{itemize}
\item if $s=0$ on a hypersurface $M$ and $Ds=0$, then $s=0$
on the whole component containing $M$.
\end{itemize}\end{re}

In the set of morphisms we can introduce an equivalence
relation:~we say that $L\sim L'$ if $L$ and $L'$ are images of
embeddings $i,i':H\hookrightarrow H_1\oplus H_2$ of a Hilbert
space $H$, such that $i-i'$ is a compact operator. If $L\sim L'$,
then $Ind_{S_1,S_2}(L)=Ind_{S_1,S_2}(L')$. If $L$ is a bordism,
then $L$ is equivalent to a direct sum of subspaces in
coordinates: $L\sim L_1\oplus L_2$, $L_i\subset H_i$, such that
$L_1$ is a finite dimensional perturbation of $\Ho_1$ and $L_2$ is
a finite dimensional perturbation of $\Hi_2$. Then
$Ind_{S_1,S_2}(L)=Ind(\Hi_1,L_1)+Ind(L_2,\Ho_2)$.

Suppose, as in Example \ref{geom-chain}, we have an elliptic
operator on a closed manifold $X$ which is composed of geometric
bordisms. Fix $n\in\Nt$ and a sequence of matrices
$$A_i\in GL_n(\B_i)\,.$$
Define a bundle $\Theta^{\{A_i\}}_X$ obtained from trivial ones on
$X_i$'s and twisted along $M_i$'s. Define bordisms $L_i(D)\in
Mor_\C(H_i,H_{i+1})$ as in Example \ref{bor-geom}.

\begin{th}\label{chain-ind} Suppose that \ref{uep}
holds for $D$ and
$D^*$ on each $X_i$ for $i=0,\dots,n$. Then
$$ind(D\otimes\Theta_X^{\{A_i\}})=n\left(\sum_{i=0}^m
Ind_{S_i,S_{i+1}}(L_i(D))\right)+
\sum_{i=1}^m \widetilde{ind}(\widetilde A_i) \,.$$
\end{th}
Here, as it was denoted before, $\widetilde A:H^{\oplus n}\to
H^{\oplus n}$ is the operator associated to the matrix $A\in
GL_n(B)$. This Theorem is a special case of Theorem
\ref{add+twist} proved below.

Taking into account Remark \ref{split-dep} the difference between
the idices of the original and twisted operator can be expressed
through the pairing  in $K$-theory.

\begin{th}\label{chain-zap}
$$ind(D\otimes\Theta_X^{\{A_i\}})-
n\,ind(D)=\sum_{i=1}^m\widetilde{ind}(\widetilde A_i) =
\sum_{i=1}^m\langle[A_i],[S_{\Hi_i}]\rangle\,.$$
\end{th}
The braked is the pairing between $[A_i]\in K^ 1(M_i)$ and
$[S_{\Hi_i}]\in K_1(M_i)$.

\section{Indices of compositions}

In \ref{split-dep} we have made some remarks about the dependence
of indices on the particular splitting. Now let us see how indices
behave under compositions of correspondences. From the
considerations in \S\ref{ioac} it is easy to deduce:

\begin{pr}\label{comp} For the composition
$$H_1\mpr{\phi}H_1\mpr{L}H_2\,,$$
where $\phi$ is a twist and $L$ is a morphism we have
$$Ind_{S_1,S_2}(L\circ\phi)=Ind_{S_1,S_2}(L)+\widetilde{ind}(\phi)\,.$$
The same holds for the opposite type composition
$$H_1\mpr{L}H_2\mpr{\phi}H_2\,,$$
$$Ind_{S_1,S_2}(\phi\circ
L)=\widetilde{ind}(\phi)+Ind_{S_1,S_2}(L)\,.$$\end{pr} On the
other hand $Ind_{S_0,S_2}(L_2\circ L_1)$ differs from
$Ind_{S_0,S_1}(L_1)+Ind_{S_1,S_2}(L_2)$ in general. This is clear
due to the basic example that comes from a decomposition
$X=X_-\cup_MX_+$. The space $L_1=H_-(D)$ is a correspondence $0\to
L^2(M;\xi)$ and $L_2=H_+(D)$ a correspondence $L^2(M;\xi)\to 0$. By 
\ref{chain-ind} we have
$$Ind_{Id,S_1}(L_1)+Ind_{S_1,Id}(L_2)=ind(D)\,,$$
while $L_2\circ L_1:0\rightarrow 0$ and $Ind_{Id,Id}(L_2\circ L_1)=0$.

Instead we have the following interesting property of indices:

\begin{th} \label{strange} The difference
$$\delta(L_1,L_2)=Ind_{S_0,S_1}(L_1)+Ind_{S_1,S_2}(L_1)-Ind_{S_0,S_2}(L_2\circ
L_1)$$
does not depend on the particular splittings.\end{th}

\Proof Since
 $$Ind_{S_{i-1},S_i}(L_i)=ind(\Hi_{i-1}\oplus
L_i\oplus \Ho_i \rightarrow H_{i-1}\oplus H_i)$$
 we have to compare indices of the operators
 $$\alpha:\Hi_0\oplus
L_1\oplus \Ho_1\oplus \Hi_1 \oplus L_2\oplus \Ho_2 \rightarrow
H_0\oplus H_1\oplus H_1 \oplus H_2$$
 and
 $$\beta:\Hi_0 \oplus L_2\circ L_1\oplus \Ho_2 \rightarrow H_0\oplus
H_2\,.$$
 The kernel  of $\alpha$
 is isomorphic to the kernel of the operator which is induced by inclusions
 $$\Hi_0\oplus L_1 \oplus L_2\oplus \Ho_2
 \rightarrow H_0\oplus H_1 \oplus H_2\,.$$
 The former operator factors through
 $$\Hi_0\oplus (L_1 + L_2)\oplus \Ho_2
 \rightarrow H_0\oplus H_1 \oplus H_2\,.$$
 Here the direct sum is replaced by the algebraic sum inside
 $H_0\oplus H_1\oplus H_2$.
 The difference of the dimensions of the
 kernels is equal to the dimension of the intersection
 $$(L_1\oplus 0)\cap (0\oplus L_2)\subset H_0\oplus H_1\oplus H_2$$
 Now we observe that the kernel of the last operator is isomorphic
 to
 $$\Hi_0 \oplus L_2\circ L_1\oplus \Ho_2 \rightarrow H_0\oplus
 H_2\,.$$
 Therefore the difference of the dimensions of the kernels of $\alpha$ and
 $\beta$ is equal to $\dim((L_1\oplus 0)\cap (0\oplus L_2))$, hence it 
does not depend on the splittings.
 We have the dual formula for cokernels and $L_i^\perp$, also not
 depending on the splittings.
\hfill$\Box$

We obtain a procedure of computing the sum of indices
$$\sum_{i=0}^m Ind_{S_i,S_{i+1}}(L_i)$$ which would not involve
splittings. We choose a pair of consecutive morphisms $L_i$,
$L_{i+1}$ and replace them by their compositions. The composition
produces a number $\delta(L_i,L_{i+1})$ and the sequence of
morphisms is shorter:
$$(L_0,L_1,\dots,L_m)\leadsto (L_0,L_1,\dots,L_i\circ
L_{i+1},\dots,L_m)+ \delta(L_i,L_{i+1})\,.$$ We pick another
composition and add its contribution to the previous one. We
continue until we get $0\rightarrow 0$. The sum of the
contributions does not depend on the splittings. One can perform
compositions in various ways. The sum of contributions stays the
same.

\begin{ex} \rm If $D$ and $D^*$ on $X_i$ and $X_{i+1}$ have the unique
extension property \ref{uep}, then $\delta(L_i,L_{i+1})=0$ as long
the gluing process along $M_{i+1}$ does not create a closed
component of $X$. If it does then $\delta(L_i,L_{i+1})$ equals to
the index of $D$ restricted to this component.
\end{ex}

\section{\label{weird}Weird decompositions of manifolds}
Let  $\{M_e\}_{e\in E}$ be  a configuration of disjoined
hypersurfaces in a manifold $X$. We assume that orientations of
the normal bundles are fixed. For simplicity assume that $X$ and
$M_e$'s are connected. Let
$$X\setminus\bigsqcup_{e\in E}M_e=\bigsqcup_{v\in V}X_v$$
be the decomposition of $X$ into connected components.
Our situation is well described by an oriented graph
\begin{itemize}
\item the vertices (corresponding to open domains in $X$) are
labelled by the set $V$
\item the edges (corresponding to hypersurfaces) are labelled by
$E$. The edge $e$ starts at the vertex $v=s(e)$ corresponding to
$X_v$ which is on the negative side of $M_e$. It ends at
$v'=t(e)$, such that $X_{v'}$ lies on the positive side of $M_e$.
The functions $s,t:E\rightarrow V$ are the {\it source} and {\it
target} functions.\end{itemize} For example the configuration

\centerline{\psfig{file=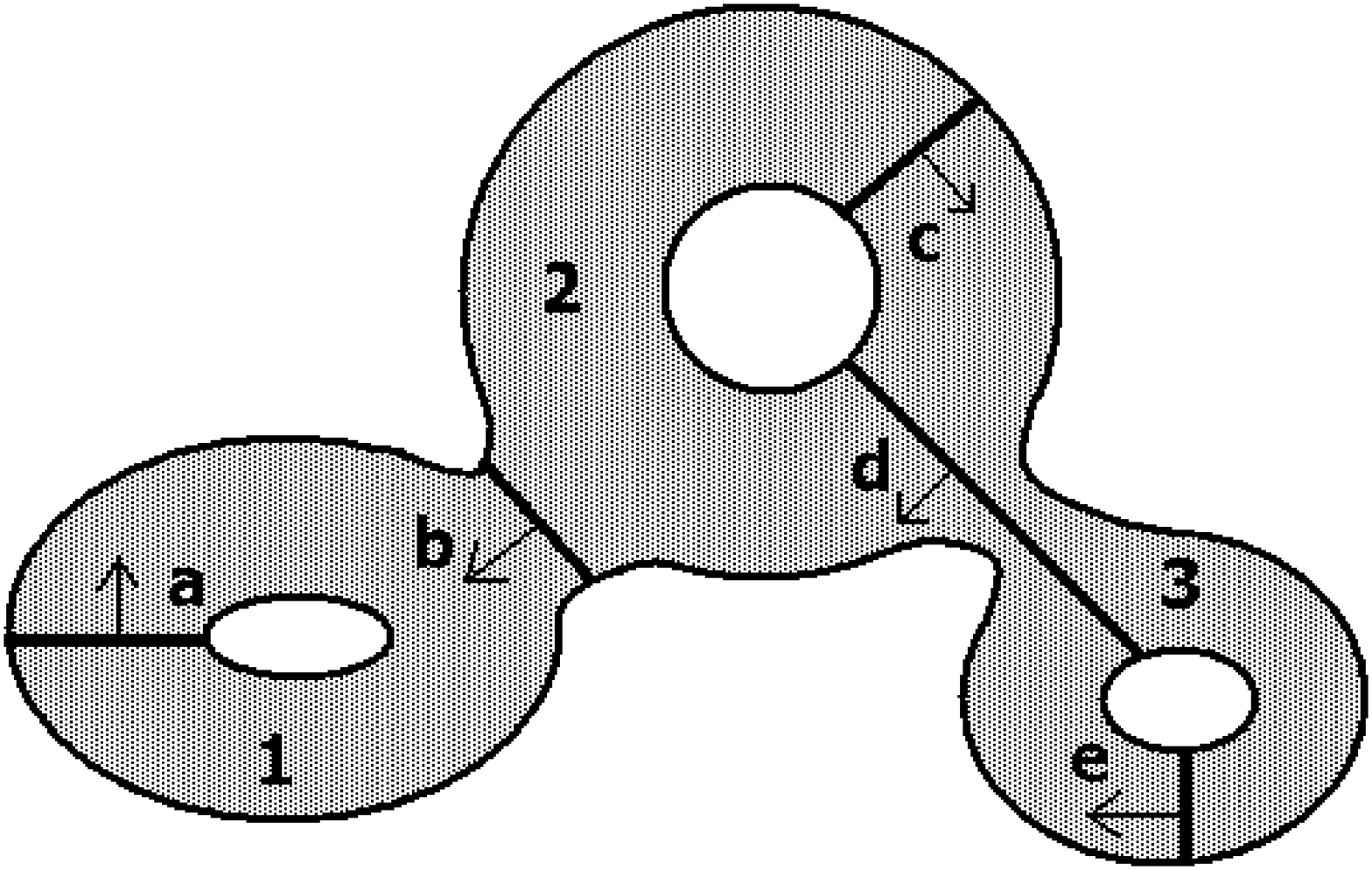,height=2truein}} \centerline{Fig.
1}

\noindent is described by the following graph:

\centerline{\psfig{file=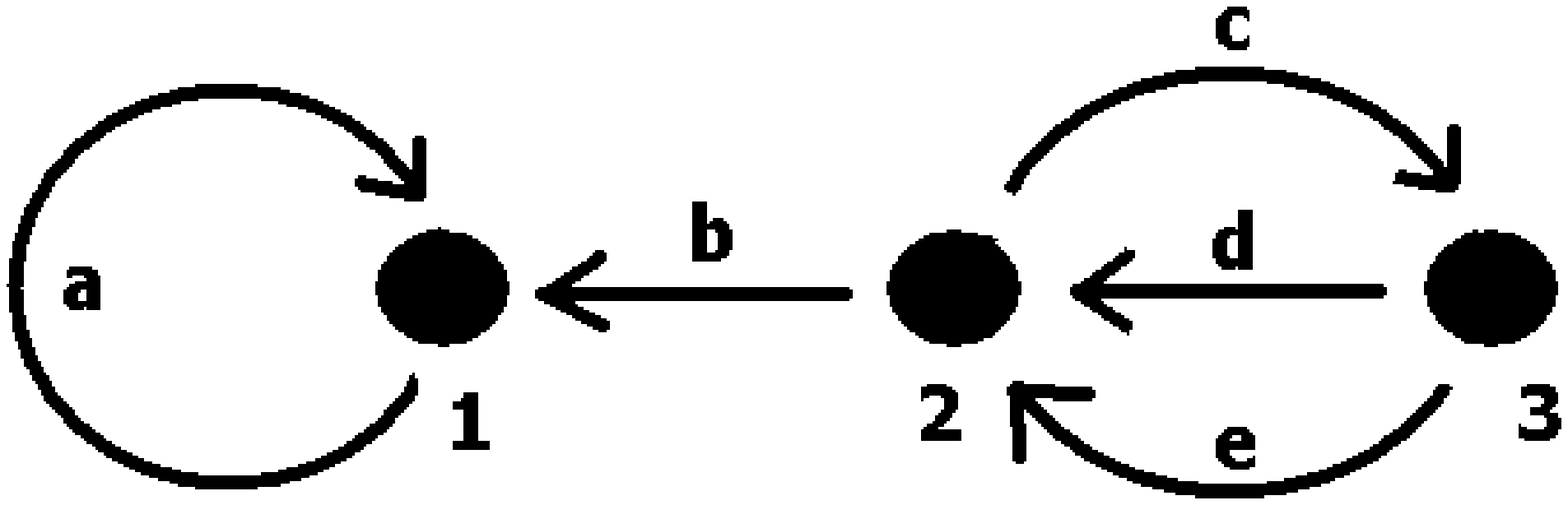,height=1.5truein}}
\centerline{Fig. 2}


\noindent A sequence of bordisms leads to the linear graph
$$\bullet_{X_0}\mpr{M_1}\bullet_{X_1}\mpr{M_2}\dots
\mpr{M_{n-1}}\bullet_{X_{n-1}}\mpr{M_{n}}\bullet_{X_n}\,.$$ Note
that this is a dual description with respect to the one presented
in Example \ref{geom-chain}. Suppose there is given an elliptic
operator $D:C^\infty(X;\xi)\rightarrow C^\infty(X;\eta)$ and a set
of transmission data $\{\phi_e\}_{e\in E}$, that is for each
hypersurface $M_e$ we are given a matrix-valued function $M_e\to
GL_n(\Ct)$.  The Riemann-Hilbert problem gives rise to the
operator
 $$D^{[\phi]}:
\bigoplus_{v\in V}C^\infty(X_v;\xi)^n\rightarrow \bigoplus_{v\in
V}C^\infty(X_v;\eta)^n\oplus \bigoplus_{e\in E}
C^\infty(M_e;\xi)^n$$
$$D^{[\phi]}(f_v)\stackrel{\rm def}{=}\left(Df_v,
\sum_{e:\, t(e)=v}f_{v|M_e}- \sum_{e:\,
s(e)=v}\phi_e(f_{v|M_e})\right)\,,\quad{\rm for}\quad f_v\in
C^\infty(X_v;\xi)^n\,.$$

For $e\in E$ let us set $H(e)=L^2(M_e;\xi)$.  The symbol of $D$
together with the choice of orientations of the normal bundles
define polarizations of $H(e)$. Let us fix particular splittings
of the spaces $H(e)$ encoded in the symetries $S_e$. Set
\begin{eqnarray*}H^{\rm bd}(v)&=&
\bigoplus_{e:\,s(e)=v}H(e)\;\oplus\bigoplus_{e:\,t(e)=v}H(e)\,,\\
H^{\rm in}(v)&=&
\bigoplus_{e:\,s(e)=v}\Ho(e)\oplus\bigoplus_{e:\,t(e)=v}\Hi(e)\,,\\
H^{\rm out}(v)&=&
\bigoplus_{e:\,s(e)=v}\Hi(e)\oplus\bigoplus_{e:\,t(e)=v}\Ho(e)\,.
\end{eqnarray*}
Let $L(v)\subset H^{\rm bd}(v)$ be the space of boundary values of
solutions of $Df_v=0$ on $X_v$. It is a perturbation of $H^{\rm
in}(v)$. For each vertex $v$ (i.e. for each open domain $X_v$) the
pair of subspaces
$$L(v),H^{\rm out}(v)\subset H^{\rm bd}(v)\,,$$
is Fredholm. Let $Ind_v$ be its index  with respect to the
polarizations $S_e$. Moreover, let
$Ind_e=Ind_{S_e,S_e}(\phi_e)=\widetilde{ind}(\phi_e)$ denote the
index of $\phi_e$, see Theorem \ref{bo}.

\begin{th}\label{add+twist}
Assume that $D$ and $D^*$ have unique extension property (\ref{uep}) on each
$X_v$. Then
 $$ind(D^{[\phi]})=\sum_{v\in V} Ind_v+\sum_{e\in E}
Ind_e \,.$$\end{th}
 In particular:

\begin{cor}\label{add} If there are no twists, i.e. each
$\phi_e=1\in GL_1(C^\infty(M_e))$, then
 $$ind(D)=\sum_{v\in V} Ind_v\,.$$\end{cor}

\noindent {\it Proof of \ref{add+twist}.} The general result
follows from the case when we have one vertex and one edge
starting and ending in it. We just sum up all $X_v$'s and all
$M_e$'s. Say that $X$ is obtained from $\hat X$ with $\partial\hat
X=M_s\sqcup M_t$ by identification $M_s$ with $M_t$.

\centerline{\psfig{file=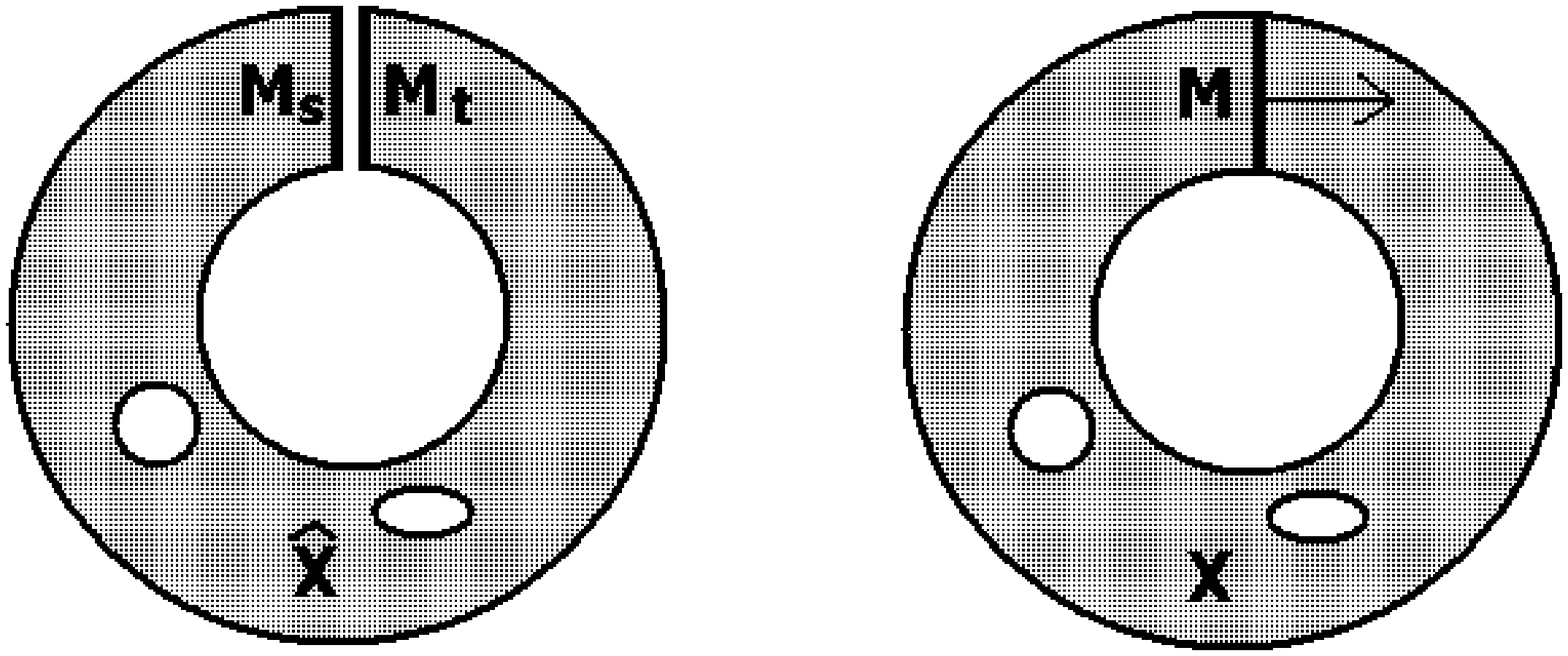,height=2truein}}
 \centerline{Fig. 3}

\noindent Then our operator $D^{[\phi]}$ is of the form:
 $$D^{[\phi]}:
C^\infty(\hat X;\xi)^n\rightarrow  C^\infty(\hat X;\eta)^n\oplus
C^\infty(M;\xi)^n$$
$$D^{[\phi]}(u)=\left(Du,
u_{|M_t}-\phi(u_{|M_s})\right)\,.$$
 We replace $\xi^{\oplus n}$ by
$\xi$ and treat $\phi$ as an automorphism of $\xi$. The index of
the operator is equal to the index of a Fredholm pair:
\begin{th} \label{heu} Let $L\subset L^2(M_s\sqcup M_t;\xi)=L^2(M;\xi)\times
L^2(M;\xi)$ be the space of boundary values of the operator $D$ on
$\hat X$. Then
$$ind(D^{[\phi]})=Ind(L,graph(\phi))\,.$$
\end{th}

The proof of Theorem \ref{add+twist} relies on this formula. We
will give a heuristic proof of \ref{heu}. The precise argument
demands introduction and consecutive use of the whole scale of
Sobolev spaces with all usual technicalities involved. The reader
may also take this formula as the definition of the index of the
problem considered above. We calculate the kernel and cokernel of
$D^{[\phi]}$:
\begin{itemize}
\item the kernel consist of solutions of $Du=0$ on $\hat X$ satisfying
$\phi(u_{|M_s})=u_{|M_t}$. By our assumption $u$ is determined by
its boundary value. Thus $$ker\,D^{[\phi]}\simeq L\cap
graph\,\phi\,.$$
\end{itemize}
The cokernel consists of
 \vskip6pt \noindent $\left\{(v,w) \in
C^\infty(\hat X;\eta^*)\oplus C^\infty(M;\xi^*) :\right.\,$\hfill

\hfill $\left.\forall u\in C^\infty(X_+;\xi)\quad\langle Du,v
\rangle + \langle u_{|M_t}-\phi(u_{|M_s}),w\rangle=0\right\}\,.$

\vskip6pt\noindent Let $G:\xi_{|M}\rightarrow\eta_{|M}$ be the
isomorphism of the bundles defined by the symbol of $D$ as in
\cite{PS}. It follows that
\begin{itemize}
\item $D^*v=0$ (since we can take any $u$ with support in $int\,\hat X$)
\item by Green formula
$\langle Du,v \rangle=\langle Gu_{|M_s},v_{|M_s}\rangle +\langle
Gu_{|M_t},v_{|M_t}\rangle$
\item since $u_{|M_s}$ and $u_{|M_t}$ may be arbitrary it follows that

 $G^*(v_{|M_s})=-\phi^*w$,

 $G^*(v_{|M_t})=w$,

\item therefore $v_{|M_s}=-G^{*-1}\phi^*G^*(v_{|M_t})$.
\end{itemize}
Now we use the identification
 $$G^*\times G^*:L^2(M_s;\eta^*)\times L^2(M_t;\eta^*)\to
 L^2(M_s;\xi^*)\times L^2(M_t;\xi^*)
 $$
under which $L^\perp$ is equal to the space of boundary values
$H(D^*)$ and
$$(graph\,\phi)^\perp=(graph(-G^{*-1}\phi^*G^*))^{op}\,.$$
(Here the opposite correspondence $R^{op}$ is defined by $(x,y)\in
R^{op}\equiv (y,x)\in R$.) In another words $\phi$ and
$G^{*-1}\phi^*G^*$ are adjoined. Since the boundary values of $v$
determine $v$ we can identify
$$coker\,D^{[\phi]}\simeq H(D^*)\cap (-graph(G^{*-1}\phi^*G^*))^{op}\simeq
L^\perp\cap(graph\,\phi)^\perp\,.$$
 \hfill$\Box$
\vskip6pt

\noindent {\it Proof of \ref{add+twist}, continuation.} After
fixing a splitting of $L^2(M;\xi)=H_e$, we have in our notation
$H^{\rm in}_v=\Hi\oplus\Ho$, $H^{\rm out}_v=\Ho\oplus\Hi$. By
\ref{bo} there exists an linear isomorphism $\Psi:H\oplus H
\rightarrow H\oplus H$ almost commuting with $\Pii\oplus\Po$, such
that $L=\Psi(\Hi\oplus\Ho)$. We parameterize the graph of $\phi$
by $\Ho\oplus\Hi$ using the composition
$\Phi=\left(\matrix{1&0\cr\phi&1\cr}\right)\circ
\left(\matrix{\Po&\Pii\cr\Pii&\Po\cr}\right)$. Thus
$$Ind(graph\,\phi,L)=ind
\left(\Phi\circ\left(\matrix{\Po&0\cr 0&\Pii\cr}\right)+
\Psi\circ\left(\matrix{\Pii&0\cr 0&\Po\cr}\right)\right)\,.$$
Since $\Psi$ almost commutes with $\Pii\oplus\Po$, the considered
operator is almost equal to the composition
$$\left(\Phi\circ\left(\matrix{\Po&0\cr 0&\Pii\cr}\right)+
\left(\matrix{\Pii&0\cr 0&\Po\cr}\right)\right)\circ
\left(\left(\matrix{\Po&0\cr 0&\Pii\cr}\right)+
\Psi\circ\left(\matrix{\Pii&0\cr 0&\Po\cr}\right)\right)\,.$$ Now
we use additivity of indices. The index of the second term is
equal to $Ind_v$. It remains to compute the first index, that is
$ind\left(\matrix{1&\Pii\cr\phi\Po&\phi\Pii+\Po}\right)$. If we
conjugate the above matrix by the symmetry
$\left(\matrix{\Po&\Pii\cr\Pii&\Po\cr}\right)$ we obtain
$\left(\matrix{\Po+\Pii\phi&0\cr\Pii+\Po\phi&1\cr}\right)$. Its
index is equal to $ind(\Po+\Pii\phi)=Ind_e$.\hfill$\Box$ \vskip6pt

The additivity of the index is not a surprise due to the well
known integral formula for the analytic index. What is interesting
in Theorem \ref{add} is that the contribution coming from separate
pieces of $X$ is also an integer number. This partition into local
indices depends only on the choice of splittings along
hypersurfaces.

\section{\label{k2}Index of a fan}

We will give another formula for the index of $D^{[\phi] }$ which
is expressed in terms of the twisted fan $\{L(i)\}$. The general
reference for fans is \cite{Bo2}. Let us first say what we mean by
a fan: it is a collection of spaces
$$L_1,L_2,\dots,L_n\subset H$$
which is obtained from a direct sum decomposition
$$H_1\oplus H_2\oplus\dots\oplus H_n=H$$ by a sequence of
twists $\Psi_1,\Psi_2,\dots,\Psi_n\in GL(H)$, i.e.
$L_i=\Psi_i(H_i)$. We assume that each $\Psi_i$ almost commutes
with each projection $P_j$ of the direct sum. We say that the fan
$\{L(i)\}$ is a perturbation of the direct sum decomposition
$H=\oplus H_i$.

\begin{th}[Index of a Fredholm fan]\label{iofp}
Let $L_1,L_2,\dots,L_n\subset H$ be a
fan. Then the following numbers are
equal:
\begin{enumerate}
\item the index of the map $\iota:L_1\oplus L_2\oplus\dots\oplus
L_n\rightarrow
H$, which is the sum of inclusions,
\item the index of the operator
$\Psi_1P_1+\Psi_2P_2+\dots+\Psi_nP_n:H\rightarrow H$,
\item the sum
$$\sum_{i=1}^n ind(P_i\Psi_i:H_i\rightarrow H_i)
=\sum_{i=1}^n ind(P_i:L_i\rightarrow H_i)\,,$$
\item the difference
$$\sum_{i=1}^{n-1}\dim(L_1+\dots+L_i)\cap L_{i+1}-\codim(L_1+\dots+L_n)\,.$$
\end{enumerate}\end{th}

\Proof The equality (1.=2.) follows from the fact that
$\Psi_i:H_i\rightarrow P_i$ is a parameterization of $L_i$. The
equality (2.=3.) follows since
$$\Psi_1P_1+\Psi_2P_2+\dots+\Psi_nP_n\sim
\prod_{i=1}^n(P_1+\dots+\Psi_iP_i+\dots+P_n)\,.$$
 To prove the
equality (1.=4.) one checks that
$$\dim(ker\,\iota)=
\sum_{i=1}^{n-1}(L_1+\dots+L_i)\cap L_{i+1} \,.$$
 This is done  by
induction with respect to $n$.\hfill$\Box$

Let us assume that the graph associated to our configuration does
not contain edges starting and ending in the same vertex (e.g.~the
situation on fig.1 is not allowed). Then $H^{\rm bd}(v)$ is a
summand in $H=\bigoplus_{e\in E}H(e)$ (there are no terms $H(e)$
appearing twice). Moreover, $\{L(v)\}_{v\in V}$ is a fan in $H$
which is a perturbation of the direct sum decomposition
$$H=\bigoplus_{v\in V} H^{\rm in}(v)\,.$$
 Consider a fan, which is twisted with respect to
$\{L(v)\}_{v\in V}$. Set $(\phi\nat L)(v)=\widetilde\phi_v(L(v))$,
where $\widetilde\phi_v$ is an automorphisms of $H$:
$$\widetilde\phi_v(f)\stackrel{\rm def}{=}\left\{
\matrix{\phi_e(f)&\;{\rm if}\;f\in H(e),\,s(e)=v\,,\cr f &\;{\rm
if}\;f\in H(e),\,s(e)\neq v\,.\cr}\right.$$

\begin{th} Assume that $D$ and $D^*$ have unique extension property
(\ref{uep}) on each $X_v$. The index of $D^{[\phi]}$ is equal to
the index of the Fredholm fan $\phi\nat L$.\end{th}

\Proof Combining Theorem \ref{add+twist} with \ref{iofp}.3 it
remains to prove that for each vertex $v$
$$ind(P_v^{\rm in}:(\phi\nat L)(v)\rightarrow
H^{\rm in}(v))=Ind_v+\sum_{e\,:\,s(e)=v}Ind_e\,.$$ If there are no
twists, then the equality follows from Proposition \ref{kryt2}. In
general the proof follows from additivity of $\widetilde{ind}$,
see Theorem \ref{bo}. \hfill$\Box$


\vskip6pt\sl

 B. B.: Institute of Mathematics PAN,

ul \'Sniadeckich 8, 00-950 Warszawa, Poland

\tt bojarski@impan.gov.pl

\vskip6pt\sl

A. W.: Institute of Mathematics, Warsaw University,

ul.Banacha 2, 02-097, Warszawa, Poland

\tt aweber@mimuw.edu.pl

\end{document}